\documentclass[11pt]{amsart}
\usepackage{amsmath,amssymb,enumerate}
\usepackage{latexsym}
\usepackage{amsfonts}
\usepackage {amsbsy}
\usepackage {amsmath}
\usepackage{amssymb}

\newcommand \C{\mathbb C}
\newcommand \N{\mathbb N}

\newcommand \R{\mathbb R}
\newcommand \B{\mathbb B}

\newcommand \U{\mathbb U}
\newcommand \T{\mathbb T}

\newcommand \sub{\subset}

\numberwithin{equation}{section}
 \title{ Subsets of full measure in a generic submanifold in $\C^n$ are non-plurithin}
 \author{Azimbay Sadullaev* and Ahmed Zeriahi}
\begin{document}

\maketitle

\noindent{\bf Abstract :} In this paper we prove that if $I\subset M $ is a subset of measure $0$ in a $C^2-$smooth generic submanifold $M \subset \C^n$, then $M \setminus I$ is non-plurithin at each point of $M$ in $\C^n$. This result improves a previous result of  A. Edigarian and J. Wiegerinck who considered the case where $I$ is pluripolar set contained in a $C^1-$smooth generic submanifold $M \subset \C^n$ ([EW10]). The proof of our result is essentially different.

 \section{Introduction}
 Real $p-$planes $\Pi \sub \C^n,
\,\,  \text{dim}_{R}\Pi=p, \,\,\, p \in \N,$ which are not contained in
any proper complex subspace of $\C^n$ are important in complex
analysis and pluripotential theory. The $\C-$hull of such plane
$\Pi$ is equal to all $\C^n$ i.e. $\Pi + J \Pi = \C^n$ ($J$ is the
standard complex structure on $\C^n$) and any non empty open
subset of $\Pi$ is non pluripolar in $\C^n$. Such planes are
called generic (real) subspaces of $\C^n$.
 Correspondingly, a real smooth submanifold $M \subset \C^n$ is said to be generic if for each $z \in M$,
its real tangent space $T_z M $ is a generic subspace of $\C^n$  i.e. $T_z M + J T_z M = \C^n$. Such submanifold has real dimension $m \geq n$.

The case of minimal dimension $\text{dim} M = n$ is the most relevant. In
this case for each $z \in M$, the tangent space $T_z M$ does not
contain any complex line i.e. $ T_z M \cap J T_z M = \{0\}$ and
$M$ is said to be {\it maximal totally real}.

One of the main notions in classical potential theory related to the Dirichlet problem is the notion of thiness.
The corresponding notion in pluripotential theory can be defined as follows.
\vskip 0.3 cm
\noindent{\bf Definition.} {\it A subset $E \subset \C^n$ is plurithin at $z^0 \in \overline{E}$ if there exists a plurisubharmonic function $u$ in a neighbourhood of $z^0$ in $\C^n$ such that
$$
u(z^0) >  \mathop
{\overline \lim }\limits_{z\rightarrow {z^0}, \,z \in {M\backslash
I}} u(z).$$
}
The main result of this paper is the following.

----------------------------------------------------------------------------------------------\\
(*) Partially supported by the fundamental research of Khorezm Mamun
Academy,Grant -1-024.

\newpage
\noindent{\bf Theorem. } {\it Let $M\subset \C^n$ be a $C^2$ generic submanifold
and let $I\subset M$ be an arbitrary  subset of measure zero in $M$. Then $M\backslash
I$ is non-plurithin at any point of $M$.}
\vskip 0.3 cm
 A. Edigarian and J. Wiegerinck proved  that if $M\subset \C^n$ be a $C^1$-smooth generic submanifold and $P \subset \C^n$ is a pluripolar set, then $M \setminus P$ is non-plurithin at any point of $M$ (see [EW10]). 

By a result of the first author ([Sa76]) for $C^3$-smooth manifolds, extended to the $C^2$-smooth case by B. Coupet ([C92]), any pluripolar set $P \subset M$ is of measure zero in $M$. Therefore for $C^2$-smooth submanifolds, our result improves the result of A. Edigarian and J. Wiegerinck. 

A powerfull method for studing the {\it pluripotential} properties on
generic submanifolds (e.g. {\it pluripolarity, pluriregularity or plurithinness}) is the well known  method of attaching analytic disks (see
[B65], [Pin74], [KhCh76], [Sa76], [C92], [EW10]). Here we  will use the same method but we need to establish more precise properties of the smooth family of analytic disks attached to $M$ along a half of the unit circle.

\vskip 0.3 cm
\noindent \emph{Acknowledgments.}  The authors whould like to thank Norm Levenberg for his helpful remarks and suggestions on a preliminary version of this paper and for providing them some useful refrences.

 This work was supported in part by
the agreement between National University of Uzbekistan and
Universit\'e Paul Sabatier de Toulouse. The authors whould like to thank these Universities for their help in organizing meetings
 of their Mathematical research groups in Tashkent and Toulouse. 

\section{Bishop's construction}

Let us recall Bishop's method for constructing analytic discs attached to real submanifolds.

 Let $M \subset \C^n$ be a totally real $C^k-$submanifold ($k \geq 1$) of dimension $n$ given locally by the
following equation
\begin {equation}
M: = \left\{ {z = x + i y \in B \times \R^n :y = h\left( x
\right)} \right\},
\end {equation}
where  $ B \subset \R^n $
 is an euclidean ball of center 0 and $
h :\ B \rightarrow \R^n$ is a $C^2$ function such that
$$ h\left( 0 \right) = 0  \ \ \text{and} \ \
Dh\left( 0 \right) = 0.
$$

 Let  $v : \T \rightarrow \R^+ $ a $C^\infty $
 function on the unit circle $\T$ such that
 $$
v|_{\gamma}  = 0 \ \
\text{and} \  \ v|_{\T \setminus \gamma\}}
> 0, \leqno (\star)
$$
where $\gamma := \left\{e^{i \theta} : \theta \in [0,\pi]\right\}$.

Assume that there exist a continuous mapping $X: \T \to \R^n $ solution of the following Bishop equation

  \begin{equation}
X\left( \tau  \right) = c - \Im \left( {h \circ X + tv}
\right)\left( \tau  \right)\,,\,\tau  \in \T,
\end{equation}
where $ \left( {c,t} \right) \in Q=Q_c \times Q_t \subset \R^n
\times \R^n $ is a fixed parameter and $\Im $ is the harmonic
conjugate operator defined by the Schwarz integral formula

         \begin{equation}
\Im \left( X \right)\left( \varsigma  \right) = \frac{1} {{2\pi
}}\int\limits_T {X\left( \tau  \right)\operatorname{Im}
\frac{{e^{i\tau }  + \varsigma }} {{e^{i\tau }  - \varsigma
}}d\tau } \,\,,\,\,\,\,\varsigma  = re^{i\theta } ,
\end{equation}
normalized by the condition
$$ \Im X\left( 0 \right) = 0.
$$
    We will consider the unique harmonic extension $
X\left( \varsigma  \right) $
 of the mapping $X\left( \tau  \right)$ to the unit disk $\U.$ Then  the following holomorphic mapping
$$
\begin{array}{ll} \Phi \left( {c,t,\varsigma } \right): = X\left(
{c,t,\varsigma } \right) + i\left[ {h^{*}\left( {c,t,\varsigma }
\right) + tv\left(
\varsigma  \right)} \right] = \\
= c + i\left\{ {h^{*}\left({c,t,\varsigma } \right) + tv\left(
\varsigma \right) + i\Im \left[ {h^{*}\left( {c,t,\varsigma }
\right) + tv\left( \varsigma \right)} \right]} \right\}
\end{array}
$$
provides a  family of analytic disks  $ \Phi \left( {c,t,\cdot
} \right):\overline \U  \to {\Bbb C}^n $ {\it attached} to $M$ along the arc $\gamma=\{e^{i\theta}: \theta \in {[0,\pi]}\}$ in the following sense
$$
\forall \left( {c,t} \right) \in Q\,,\, \ \ \forall \, \tau  \in \gamma
\,,\, \ \ \Phi \left( {c,t,\tau } \right) \in M.
$$
Here  $ X\left( {c,t,\cdot} \right),\,\,\,h^{*}\left(
{c,t,\cdot} \right)\,\text{and}\,\,v $
 denote the harmonic extensions of  $
X\left( {c,t,\tau } \right)
\newline\,\, h \circ X\left( {c,t,\tau }
\right)\,\,\,\text{and}\,\,\,v\left( \tau  \right) $
from  the unit circle $\T$  to the unit disk $\U$, respectively.

  We need a smooth family $\Phi \left( {c,t,\cdot } \right)$ of analytic discs. This is provided by the following result of B.Coupet ([C10]) extending a construction done previously by the first author.
  \vskip 0.2 cm

\noindent{\bf Theorem ([C92])}. {\it  Let $h\in C^k (\B),\,\,k \in \N^*$. Then there exist a
neighborhood $Q= Q_c \times Q_t \ni 0$  such, that the Bishop
equation (2.2)  has a unique solution $u \in
C^{k-1}(Q \times \T).$}

 {\it Moreover, the harmonic extensions of $\,\,u$, $h\circ u\,\,$ to
 the unit disk $\U$, their conjugates $\Im u, \Im h\circ u$
belong  $\,\,C^k{(Q \times \U)}\cap C^{k-1}{(Q \times \overline \U)}$
and 
$$\left\| \Im u \right\| \leq A \left\| u \right\|, \left\| \Im
h \right\|\leq A \left\| h \right\|, 
$$
 where $A$ is a constant.}
\vskip 0.2 cm

Therefore, when the submanifold $M$ is $C^1-$smooth, the solution
$X(\tau,c,t)$ is continuous in $Q \times  \overline {\U}$ and for
twise smooth submanifold we obtain a $C^1-$smooth familly of disks, attached
to $M$.

In the case when $I\subset M$ is pluripolar, Edigarian and
Wiegerinck [EW] needed only a continuous family of disks and gave a
beautiful proof for the non-thinness  of $M\backslash I$.

Recall that every pluripolar set $I\subset M$  on the twice smooth
generic manifold $M$ has zero-measure (see [Sa76], [C92]), but the converse
is far from being true. Indeed, there are many subsets of $\R^n \subset \C^n$  with zero-measure which are not pluripolar in $\C^n$.
So our theorem is a non trivial improvement of Edigarian and Wiegerinck theorem in the case of $C^2-$manifolds.
We need to assume $C^2-$smoothness of $M$ because we need smoothness of the family of discs constructed in ([C92]).

   Observe that the family of discs $\Phi (c,t,\zeta)$ obtained above
 \begin {equation}
 \Phi (c,t,\zeta) = X (c,t,\zeta) + i (h^{*} (c,t,\zeta) + t v (\zeta)), \ (c,t) \in Q, \zeta \in \bar \U,
 \end {equation}
 satisfies the following conditions:

 \begin{equation} \label{eq:PFE}
  X (c,t,\tau) = c - \Im \left( h \circ X (c,t,\tau) + t v(\tau) \right),  \,\,  (c,t) \in Q, \,\,  \tau \in \partial \U.
 \end{equation}

\begin{equation}
h^{^{*}}(c,t, \tau)=h \circ X (c,t, \tau), (c,t) \in Q, \tau \in
\partial \U
\end{equation}

\begin{equation}\begin{array}{ll}
X(c,0,\zeta)\equiv c, h^{^{*}}(c,0, \zeta)\equiv h(c) \,\,\,
\text{so
that} \\
\Phi(c,0, \zeta)\equiv c+ ih(c) \in M, c \in Q_c
\end{array}
\end{equation}

\begin{equation}
X(c,t,0)= \frac{1}{2\pi} \int_{T}X(c,t, \tau) d \tau \equiv c,
\,\, (c,t) \in Q.
\end{equation}

 \section {Proof of the theorem }

 \noindent{\bf 1. First case.}  We consider the case when $\text{dim}\,M=n$. Fix a point $p$,
say $p = 0 \in M$. We want to prove that $M\backslash I$ is non-plurithin
at $0$. Assume the contrary, that there exist a $V(z)\in {PSH}(G)$ such that $\,\,V(0)> 1/2$,
$V|_{M \backslash
{I\cup\{0\}}}=0$ and $\,\,0\leq V(z) \leq 1,\,\,$ for $z \in G$,  where $G$ a neighbourhood of
$0\in \C^n$.

The proof of the theorem goes in several steps.

 ${\bf 1^0}$. First, for convenience, we introduce the following terminology. We say that the disc
  $\Phi_{c,t} := \Phi(c,t,\cdot)$, defined by (2.4) is {\it good} if
 $$
 H_1 (\gamma_I) = 0, \ \text{where} \ \gamma_I := \{\tau \in \gamma: \Phi (c,t,\tau) \in I\}.
 $$
 Assume that for some fixed value of the parameter $(c,t) \in Q $, the corresponding disc $ \Phi_{c,t}$ is
  "good".  Then the function defined on $\U$ by the formula
 $u (\zeta) := V \circ \Phi (c,t,\zeta)$ satisfies the following properties (we note that $\Phi (c,t, \cdot): \, \overline {U} \to G$ for a small enough  $\,Q$):
  $$
  u \in \,{SH}(\U)\cap C \left(\overline{\U} \cup \{\gamma \backslash \gamma_I \}\right),\,\,u|_{\gamma \backslash
   \gamma_P}=0.
   $$

 Let $ \omega (\zeta,\gamma \setminus \gamma_I,\U)$ be the harmonic measure of the set $\gamma \setminus \gamma_I$ with
  respect to $\U$ at the point $\zeta$. Then $\omega (\zeta,\gamma \setminus \gamma_I,\U)$ is the generalized solution
   of the Dirichlet problem in $\U$ with boundary data equal to $- \chi_{\gamma\backslash\gamma_{I}}$ on
   $\partial \U$, where $\chi_{\gamma\backslash\gamma_{I}}$ is the chararcteristic function of the Borel
   set $\gamma \setminus \gamma_I$.
 This means that $\omega (\zeta,\gamma \setminus \gamma_I,\U)$ is a harmonic function negative on $\U$ and
 equal to$- \chi$ quasi-everywhere on $\partial \U$, in particular  $\omega (\zeta,\gamma \setminus \gamma_I,\U) = - 1$
 on $\gamma \setminus \gamma_I$.  Since $H_1 (\gamma_I) = 0$, it follows from Poisson integral formula that
 $ \omega (\zeta,\gamma \setminus \gamma_I,\U) \equiv \omega (\zeta,\gamma ,\U)$, since they coincide almost everywhere
  on $\partial \U$.

We put
$$
\Omega  = \left\{ \zeta  \in \overline \U : 1 + \omega \left( \zeta
,\gamma ,\U \right) < 1/4 \right \} .
$$

Then $ \Omega  \supset \gamma \,\,\text{and}\,\,\Omega  = \mathop
\Omega \limits^ \sim   \cap \overline \U $
 for some open set $
\mathop \Omega \limits^ \sim  :\gamma  \subset \mathop \Omega
\limits^ \sim $. We note that $ \Omega $
 depended only of $
\gamma .$
     By the so called two-constant theorem, it follows that for any $
\varsigma  \in \U $

 \begin{equation}
0 \leqslant u\left( \zeta  \right) \leqslant  {1 + \omega \left(
{\zeta ,\gamma ,\U} \right)}  .
\end{equation}

This implies that $ u \equiv 0 $
 on $\gamma
$
 and then
 \begin{equation}
V( z) < 1/4,
\end{equation}
{\it provided} that the point $z$ lies in the image $ \Phi_{c,t} \left( \Omega  \right) $
 of the good disk  $
\Phi_{c,t}. $\\

$ {\bf 2^0}$. The next step is to prove that for arbitrary fixed
$\zeta^0  \in \Omega  \cap \U $ the images $ \Phi \left(c,t,\zeta
^0\right),\,\,$ as $(c,t)$ vary in $Q$, fill up an open set $ W(\zeta^0) \ni 0$ in $\C^n = \R^n \times \R^n$.

Consider the map
$$
S:\left( {c,t} \right) \to \Phi \left( {c,t,\zeta^0 } \right) =
X\left( {c,t,\zeta^0 } \right) + i[h^{*}\left( {c,t,\varsigma ^0 }
\right) + tv\left( {\varsigma ^0 } \right)]\,.
$$
The jacobian of S is

 $$
J(c,t)= \left| \begin{matrix}
   {D_{c_j } X} & {...} & {D_{t_j } X}  \\
   {...} & {...} & {...}  \\
   {D_{c_j } h^{*}} & {...} & {v\left( {\zeta^0 } \right) \mathbb I + D_{t_j } h^{*}}  \\
 \end{matrix}  \right|
 $$
where  $ \mathbb I$ is the identity matrix in $\R^n$. We see that
 $$
 J(c,t) = v^n\left( {\varsigma ^0 } \right)\frac{{\partial X_1 }}
 {{\partial c_1 }}...\frac{{\partial X_n }} {{\partial c_n }} +
 \,O\left(\sum_{i \neq j} \Vert D_{c_j} X_i \Vert , \sum_j \left[\Vert D_{c_j} h^*\Vert +
 \Vert D_{t_j} h^*\Vert\right]\right)
  $$
  Here and belove $ \Vert  \cdot  \Vert\ $ denote the $sup$    norms.

If  $t=0$ then by (2.7)  $ X\left( {c,0,\zeta^0 } \right)
 \equiv c $
 and
 $
 h^{*}\left( {c,0,\zeta^0 } \right) \equiv h\left( c \right)
$. Hence by smoothness of $\Phi$ and the condition $(\star)$, it follows that
$$ J \ne 0
$$
 for small enough  $Q _c $ and $Q_t$ (remember that $h(0)=Dh(0)=0).$ It follows that the
transformation $S $ is open and then the following set 
$$
W(\zeta^0)= S (Q) =  \left\{\Phi \left(c,t,\zeta^0\right) ; (c,t) \in \left( Q \right) \right\}
$$  is open in $\C^n = \R^n \times \R^n$ and contains $ 0$.

 Shrinking, if
necessary, $Q\,\,$ we can assume, that $W(\zeta)\subset
G, \,\,\forall \,\zeta \,\in \Omega \cap \U$ and put $W=\bigcup
\left \{W(\zeta ): \,\,\zeta\in {\Omega \cap \U}\right\}.$\\

${\bf 3^0}$.  Set
 $$
W' = \bigcup_{c,t} \left\{W \cap \Phi_{c,t} \left(
\Omega \right) : \, \,\,\Phi_{c,t} \, \, \, \text{is a good disc}\right\}.
$$
  Then (3.2) implies that  $ V(z) < 1/4 $ for $z \in W'$.

 Assume for the moment that $ H_{2n} \left( {W\backslash W'} \right) = 0$ (see next steps $\,\, 3^0-5^0$). Then it follows that from the submean value property of subharmonic functions in $\R^{2 n}$ that $V(z) < 1/4 $
 everywhere on $W^{\circ}$, i.e. for any $z^{0} \in W^{\circ}$,
$$
\limsup_{z \to z^0, \, \, z \in W'} V(z) = V(z^0)
$$
 In particular, $V(0)< 1/4$. This contradiction will prove the theorem in the case when $\text{dim}\, M=n$.

  The question remains to know "how much" good disks we have?
  In fact we have enough good discs to prove that $H_{2 n} (W \setminus W') = 0$. Indeed,
   we calculate  $
D_\tau  X\left( {c,t,\tau } \right),\,\,\,\,\left( {c,t} \right)
\in Q, \,\,\,\,\tau  \in T $
 \begin{equation}
D_\tau  X\left( {c,t,\tau } \right) =  - D_\tau  \Im h \circ
X\left( {c,t,\tau } \right) - tD_\tau  \Im v\left( \tau  \right)
\end{equation}

Since, $ D_\tau  \Im v\left( \tau  \right) = D_{\mathop n\limits^
\to  } v\left( \tau  \right) $, where   $ D_{\mathop n\limits^ \to
} $
  is the normal derivative $
\mathop n\limits^ \to $, then 
$$
D_\tau  X\left( {c,t,\tau } \right) + tD_{\mathop n\limits^ \to  }
v\left( \tau  \right) = \Im D_\tau  h \circ X\left( {c,t,\tau }
\right).
$$
For   $k$-coordinate of  vector $X\left( \tau  \right) = X\left(
{c,t,\tau } \right) $
 we have
 \begin{equation}\begin{array}{ll}
\left\| {D_\tau  X_k \left(c,t, \tau  \right) + t_{k}D_{\mathop
n\limits^ \to  } v\left( \tau  \right)} \right\| = \left\| {\Im
D_\tau  h_{k} \circ X\left(c,t, \tau  \right)}
\right\|  \leqslant \\
 \leqslant {const}\left\|
{D_\tau  h_{k} \circ X\left(c,t, \tau  \right)} \right\| \leqslant
O\left( \varepsilon  \right)\left\| {D_\tau X\left(c,t, \tau
\right)} \right\|,
\end{array}
\end{equation}
where $\varepsilon =\left\{\left\|c\right\|+ \left\|t\right\|:
\,\,c\in Q_c,\,t\in Q_t \right\}$ and
$\left\|(a_1,....,a_n)\right\|= {\max}\{|a_1|,...,|a_n|\}$.

Therefore,
\begin{equation}\begin{array}{ll}
\left| {t_k D_{\mathop n\limits^ \to  } v\left( \tau \right)}
\right| - O\left( \varepsilon  \right)\left\| t \right\| \leqslant
\left| {D_\tau  X_k \left( {c,t,\tau } \right)} \right| \leqslant
 \\ \leqslant \left| {t_k D_{\mathop n\limits^ \to  } v\left( \tau
\right)} \right| + O\left( \varepsilon  \right)\left\| t
\right\|\,\,,\,\,1 \leqslant k \leqslant n\,\,,\,\,\tau  \in T
\end{array}
\end{equation}

The second part of (3.5) implies
 \begin{equation}
\left\| {D_\tau  X\left( {c,t,\tau } \right)} \right\| \leqslant
C\left\| t \right\|\,\,,\,\,\left( {c,t,\tau } \right) \in Q
\times T,
\end{equation}
where $C > 0$ is a constant.

If we denote by $ b = \mathop {\inf }\limits_{\gamma ^0}  \left|
{D_{\mathop n\limits^ \to  } v\left( \tau  \right)} \right| > 0, $ where $\gamma ^0 \subset \subset \gamma,$
  then decreasing  $
\varepsilon  > 0 $ in (3.5)  we can arrange so that $O\left( \varepsilon  \right) < \frac{b} {2}$ and the first part of $(3.5)$ implies
\begin{equation}
 \vert D_{\tau} X_k (c,t,\tau)\vert \geq \vert t_k\vert b - \Vert t\Vert b \slash 2,
\end{equation}
for $\tau \in \gamma^0, 1 \leq k \leq n$.\\

${\bf 4^0 }$.  Fix  $ \left( {x^0 ,y^0 } \right) \in W\backslash
M$. Then by $p.\, 2^0 $
 there exist $
\left( {c^0 ,t^0, \zeta^0} \right) \in Q\times (\Omega \cap \U) $
 such that  $
\Phi \left( {c^0 ,t^0 ,\zeta^0 } \right) = \left( {x^0 ,y^0 }
\right) $.
 Let $\|t^0\|=|t^{0}_{k}|$, for simplicity we assume
$k=n$, and let $ 'c = \left( {c_1 ,...,c_{n - 1} }
\right)\,\,,\,\,'t = \left( {t_1 ,...,t_{n - 1} } \right) $.

We consider the transformation
 \begin{equation}
S:\left( {'c_{} ,'t,\varsigma } \right) \to \Phi \left( {'c_{}
,c_n^0 ,'t,t_n^0 ,\varsigma } \right): {'Q} \times \overline \U
\rightarrow \C^n,
\end{equation}
where $$ 'Q := Q \cap \{c_n = c_n^0, t_n = t_n^0\} \subset \R^{2n
- 2}.
$$
Then $ S\left( {'c^0 _{} ,'t^0 ,\varsigma ^0 } \right) = \left(
{x^0 ,y^0 } \right) $
  and its jacobian is equal to

 $$
J('c,'t,\zeta) = \left| {\begin{matrix}
   {\frac{{\partial x_1 }}
{{\partial c_1 }}} & {...} & {\frac{{\partial x_{n - 1} }}
{{\partial c_1 }}} & | & {\frac{{\partial y_1 }} {{\partial c_1
}}} & {...} & {\frac{{\partial y_{n - 1} }} {{\partial c_1 }}} & |
& {\frac{{\partial x_n }} {{\partial c_1 }}} & {\frac{{\partial
y_n }}
{{\partial c_1 }}}  \\
   {...} & {...} & {...} & | & {...} & {...} & {...} & | & {...} & {...}  \\
   {\frac{{\partial x_1 }}
{{\partial c_{n - 1} }}} & {...} & {\frac{{\partial x_{n - 1} }}
{{\partial c_{n - 1} }}} & | & {\frac{{\partial y_1 }} {{\partial
c_{n - 1} }}} & {...} & {\frac{{\partial y_{n - 1} }} {{\partial
c_{n - 1} }}} & | & {\frac{{\partial x_n }} {{\partial c_{n - 1}
}}} & {\frac{{\partial y_n }}
{{\partial c_{n - 1} }}}  \\
   { -  - } & { -  - } & { -  - } & | & { -  - } & { -  - } & { -  - } & { -  - } & { -  - } & { -  - }  \\
   {\frac{{\partial x_1 }}
{{\partial t_1 }}} & {...} & {\frac{{\partial x_{n - 1} }}
{{\partial t_1 }}} & | & {\frac{{\partial y_1 }} {{\partial t_1
}}} & {...} & {\frac{{\partial y_{n - 1} }} {{\partial t_1 }}} & |
& {\frac{{\partial x_n }} {{\partial t_1 }}} & {\frac{{\partial
y_n }}
{{\partial t_1 }}}  \\
   {...} & {...} & {...} & | & {...} & {...} & {...} & | & {...} & {...}  \\
   {\frac{{\partial x_1 }}
{{\partial t_{n - 1} }}} & {...} & {\frac{{\partial x_{n - 1} }}
{{\partial t_{n - 1} }}} & | & {\frac{{\partial y_1 }} {{\partial
t_{n - 1} }}} & {...} & {\frac{{\partial y_{n - 1} }} {{\partial
t_{n - 1} }}} & | & {\frac{{\partial x_n }} {{\partial t_{n - 1}
}}} & {\frac{{\partial y_n }}
{{\partial t_{n - 1} }}}  \\
   { -  - } & { -  - } & { -  - } & | & { -  - } & { -  - } & { -  - } & | & { -  - } & { -  - }  \\
   {\frac{{\partial x_1 }}
{{\partial \varsigma '}}} & {...} & {\frac{{\partial x_{n - 1} }}
{{\partial \varsigma '}}} & | & {\frac{{\partial y_1 }} {{\partial
\varsigma '}}} & {...} & {\frac{{\partial y_{n - 1} }} {{\partial
\varsigma '}}} & | & {\frac{{\partial x_n }} {{\partial \varsigma
'}}} & {\frac{{\partial y_n }}
{{\partial \varsigma '}}}  \\
   {\frac{{\partial x_1 }}
{{\partial \varsigma ''}}} & {...} & {\frac{{\partial x_{n - 1} }}
{{\partial \varsigma ''}}} & | & {\frac{{\partial y_1 }}
{{\partial \varsigma ''}}} & {...} & {\frac{{\partial y_{n - 1} }}
{{\partial \varsigma ''}}} & | & {\frac{{\partial x_n }}
{{\partial \varsigma ''}}} & {\frac{{\partial y_n }}
{{\partial \varsigma ''}}}  \\

 \end{matrix} } \right|,
$$

Here $\varsigma = \varsigma' + i \varsigma''$ and
$$
 x_k ('c,'t,\zeta) = X_{k}\left({'c_{} , c_n^0,'t,t_n^0,\zeta}\right),  \ \ k = 1, \cdots n,
$$
$$
 y_k ('c,'t,\zeta) = h_k^* \circ X\left( 'c_{} , c_n^0,'t,t_n^0,\zeta\right) + t_k v\left(
\zeta\right), \ \ k = 1, \cdots n - 1,
$$
$$
y_n ('c,'t,\zeta) = h_n^*  \circ X\left( 'c_{},
c_n^0,'t,t_n^0,\zeta \right) + t_n^0 v\left( \zeta \right)).
$$

The determinant $J$,  is composed by $9$ block matrices $D_{ij}
\,\,,\,\,i,j = 1,2,3$.

We will show that $ J ('c^0,'t^0,\zeta^0) \neq 0$, which will
imply that  the operator  $S $
 is local diffeomorphism in a neighbourhood of the point $\left( {'c^0 ,'t^0 ,\zeta ^0 }\right)\,$

   By (2.7)  $
X\left(c,0,\zeta\right) \equiv c\,\,,\,\,h^{*}\left( {c,0,\zeta}
\right) \equiv h\left( c \right) $. Therefore,

$$
 \left| {\begin{matrix}
  {D_{1 1}} &  {D_{1 2}}  \\
   {D_{2 1}} & D_{2 2}  \\
  \end{matrix} } \right|_{('c,0,\zeta)} =  D_{11} .D_{22}  = v^{n - 1} \left( {\zeta } \right)
$$
and
$$
\left| {\begin{matrix}
  {D_{1 1}} &  {D_{1 2}}  \\
   {D_{2 1}} & D_{2 2}  \\
  \end{matrix} } \right|_{('c,'t,\zeta)} =  v^{n - 1} \left( {\zeta } \right) +
  O\left(\varepsilon\right).$$
 Note also that if $\zeta = \zeta' + i \zeta''$ then

\begin{equation}
D_{3 3} =  \left| {\begin{matrix}
  {\frac{{\partial x_n }}{{\partial \zeta' }}} &  {\frac{{\partial y_n }}
{{\partial \zeta'}}}  \\
   {\frac{{\partial x_n }}
{{\partial \zeta''}}} & {\frac{{\partial y_n }}
{{\partial \zeta'' }}}  \\
  \end{matrix} } \right| = \left|\frac{ d }{d \zeta}(x_n + i y_n)\right|^2.
\end{equation}
 Now consider the right hand side near the arc$\gamma^0 $. It is clear that for every $s > 0,$ there is an
 open set $\Omega' \supset \gamma^0$ such that
 $$
 \left|\frac{ d }{d \zeta}(x_n + i y_n) ('c,'t,\zeta) \right|^2 \geq \left| D_{\tau} x_n ('c,'t,\tau)\right|^2 - s,
 \forall \zeta \in \U \cap \Omega', \tau \in \gamma^0, ('c,'t) \in Q.
 $$

 By $(3.6)$, $(3.7)$ and $(3.9)$,  it follows that
 \begin{eqnarray*}
 \left|J ('c,'t,\zeta)\right| &= & \left|D_{1 1}|\cdot |D_{2 2}\right| \cdot \left|\frac{ d }{d \zeta}(x_n + i y_n)
  ('c,'t,\zeta)\right|^2 + O (\varepsilon) \\
 &\geq & \left[v^{n - 1} (\zeta) + O (\varepsilon)\right] \cdot \left[|t_n b \slash 2|^2
   - s \right] + O (\varepsilon),
 \end{eqnarray*}
 for all $('c,'t,\zeta) \in {'Q} \times \left[\U \cap \Omega'\right]$.

 We can take $\Omega \cap \Omega'$ instead of $\Omega$ and observe that all functions $O (\varepsilon)$ do not depend
 on $\xi.$
 Therefore if we take $\varepsilon,  s$ small enough, then $ \left|J ('c^ 0,'t^ 0,\xi^ 0)\right| > 0$ and, in particular,
 the operator
 $$
 S ('c,'t,\zeta) : {'\hat Q} \times \{\vert \zeta - \zeta_0\vert < \sigma'\} \rightarrow \U (x^ 0,y^ 0)
 $$
is an homeomorphism, where $\sigma' > 0$ and ${'\hat{Q}} =
{'\hat{Q}}_c \times  {'\hat{Q}}_t \subset {'Q}$ is a neighbourhood of
$('c^ 0,'t^ 0)$ and $U (x^0,y^0)$ is a neighbourhood of
$(x^0,y^0)$.

 ${\bf 5^0.}$ We show that there exist a neighborhood  $\check{Q}_t \subset\,{'\hat{Q}}_t$ such that for every fixed
   $'t\in {\check{Q}_t}$ the mapping
 $$
    S\left( {'c,'t ,\tau } \right)=\Phi \left('c,c^0_n,'t,t^0_n, \tau \right):\,\,'Q_c \times \gamma^0 \to
M$$ is local homeomorphism. We put

$$
A_k  = \left( {\begin{matrix}   {\frac{{\partial x_1 }}
{{\partial c_1 }}} & {...} & {\frac{{\partial x_n }}
{{\partial c_1 }}}  \\
   {...} & {...} & {...}  \\
   {\frac{{\partial x_1 }}
{{\partial c_{n - 1} }}} & {...} & {\frac{{\partial x_1 }}
{{\partial c_{n - 1} }}}  \\
 \end{matrix} } \right)^{\surd k},
$$
where $ \left( {...} \right)^{\surd k} $ means that the
$k^{\text{th}}$ column is omitted.  Then $ A_1  = ...=A_{k - 1}  =
0\,\,,\,\,A_n  = 1\,\,\text{if}\,\,t = 0 $. We consider the
jakobian-minor
$$
J('c, 't, \tau) = \bmod \left| {\begin{matrix}
   {\frac{{\partial x_1 }}
{{\partial c_1 }}} & {...} & {\frac{{\partial x_n }}
{{\partial c_1 }}}  \\
   {...} & {...} & {...}  \\
   {\frac{{\partial x_1 }}
{{\partial c_{n - 1} }}} & {...} & {\frac{{\partial x_n }}
{{\partial c_{n - 1} }}}  \\
   {\frac{{\partial x_1 }}
{{\partial \tau }}} & {...} & {\frac{{\partial x_n }}
{{\partial \tau }}}  \\
 \end{matrix} } \right| \geqslant \left| {\frac{{\partial x_n }}
{{\partial \tau }}} \right|.\left| {A_n } \right| - \left( {\left|
{\frac{{\partial x_1 }} {{\partial \tau }}} \right|\left| {A_1 }
\right| + ... + \left| {\frac{{\partial x_{n - 1} }} {{\partial
\tau }}} \right| \left| {A_{n - 1} } \right|} \right)$$ $$=
\left| {A_n } \right| \left\{  \left| {\frac{{\partial x_n }}
{{\partial \tau }}} \right|  - \left( { \left| {\frac{{\partial
x_1 }} {{\partial \tau }}} \right|} {\frac{{ \left| {A_1 }
\right|}} {{\left| {A_n } \right| }}} + ... +{\left|
{\frac{{\partial x_{n-1} }} {{\partial \tau }}} \right|} {\frac{{
\left| {A_{n-1} } \right|}} {{\left| {A_n } \right| }}}
 \right) \right\}.
$$
By (3.6) and   (3.7) we have
$$\left\| {D_\tau  X\left( {c,t,\tau } \right)} \right\| \leqslant
C\left\| t \right\|\,\,,\,\,\left( {c,t,\tau } \right) \in Q
\times T\,\,,\,\,C - \text{constant},$$
  $$ \vert D_{\tau} X_n (c,t^0,\tau)\vert \geq \vert t^0_n\vert b - \Vert t^0\Vert b \slash 2=
                       \left( {b - \frac{b} {2}} \right)\left| {t_n^0 } \right| =
\frac{{b\left\| {t^0 } \right\|}} {2},\,\,\left( {c,t,\tau }
\right) \in Q \times T.$$
  Therefore,   $ J('c,'t^0,\tau)\neq 0 $ if $t^0$ small. The same is true
for a neighborhood $\check{Q}_t \subset\,'\hat{Q}_t.$

    Since the set $I\subset M$ has zero-measure, then for almost
    every $'c \in {'Q}_c$ the disks $\Phi\left('c, 't, \tau \right), \text{with fixed} \,\, 't \in \check{Q}_t$, intersects $I$ on zero-lenght set, i.e. the disks
$\Phi\left('c, 't, \zeta \right)$ are "good" for almost every $'c
\in {'Q}_c.$
Therefore, the $ \left( {W\backslash W'} \right) \cap
\left\{ {\Phi \left( {'c,{'t},\varsigma } \right) \in
\mathbb{R}^{2n} :\, 'c \in {'Q_1} ,\varsigma  \in B\left( {\varsigma
^0 ,r} \right)} \right\},\,r>0 $, has zero $ \left( {n + 1}
\right)$-measure. It follows that in a neighbourhood of $ \left(
{x^0 ,y^0 } \right) \in W $ the set $W\backslash W' $
 has zero $2n$-measure, which completes the proof that $H_{2 n} (W \setminus W') =
 0$, since $(x_0,y_0)$ is arbitrarily fixed.\\

 \noindent {\bf 2. General case}. Let $M$ be an arbitrary generic manifold of dimension $m > n$ and let $I\subset M$ be a subset of measure zero in $M$.

  Fix a point, say $z^0 = 0 \in M$.
 Changing holomorphic coordinates in $\C^n$, we can assume that the tangent space $T_0 M$,
 which by definition does not contain any complex
 hyperplane, can be written as
 $$
 T_0 M = \{z = x + i y \in \C^n : y_1 = \cdots = y_{2 n - m} = 0\}.
 $$
 Hence for a small neighbourhood $G = G_1 \times G_2$ of the origin with
 $$G_1 =\{ (x,y'') = (x,y_{2 n - m + 1}, ... , y_n) \in \R^n \times \R^{m - n} : \vert x\vert \leq \delta,
 \vert y''\vert < \delta \},$$
 $$G_2 =\{ y' = (y_1, \cdots, y_{2n - m}) \in \R^{2 n - m} : \vert y'\vert < \delta\},$$

 we can represent $M$ as a graph
 $$
  M \cap G = \{z  \in G :  y' = h (x,y'') \},
 $$
 where $h$ is $C^2$ smooth  mapping from $G_1$ into $\R^{2 n - m}$.

 Observe that  for each $ (m-n) \times n$ matrix $B$ the intersection $M\cap\Pi_B$ of $M$ with the plane
 $\Pi_B := \{ z \in \C^n : y'' = x \cdot B\}$ is
 an $n-$dimensional generic manifold for any small enough $B$.
 Since $H_m (I) = 0$ then there exist $B$ such that  $I\cap \Pi_B$ has zero-measure.
 Hence there exists at least one generic submanifold $M'$ of dimension $n$ such that $\,\,0\in M' \subset M$ and  $H_n (I \cap M') = 0.$

 By the first case the set $M'\backslash I$ is not thin at $0$. It follows, that $M\backslash I$ also is not thin at $0$.$ \rhd$\\

\noindent{\bf Open problem :} Let $M \subset \C^n$ be a smooth generic submanifold of dimension $m$ and $E \subset M$ a given Borel subset such that $ 0 \in E$ and the following condition holds
$$
\lim_{r \to 0} \frac{H_m (E \cap \B_r)}{r^m} =0,
$$
where $H_m$ is the Hausdorff measure of dimension $m$ and $\B_r \subset \C^n$ is the euclidean ball of radius $r$. 
Then is-it true that the set $M \setminus E$ non-plurithin at $0$? \\

                     \begin{center}
                     \textbf{  References}

                          \end{center}
\vskip 0.5 cm

\noindent [B65] E.Bishop: {\it Differentiable manifolds in complex Euclidean spaces}. Duke Math.J., V.32 (1965), no.1, 1-21.

\noindent [CPL05] D.Coman, N.Levenberg and E.Poletsky: {\it Smooth submanifolds intersecting any analytic curve in a
    discrete set}. Math. Ann. 332 (2005), no. 1, 55-65.

\noindent[C92] B. Coupet: {\it Construction de disques analytiques et r\'egularit\'e de fonctions holomorphes au bord.} Math.Z.
 209 (1992), no.2, 179-204.

\noindent[EW10] A. Edigarian and J. Wiegerinck: {\it  Sherbina's theorem for finely holomorphic functions.} Math.Z. 266 (2010), no.2, 393-398.

\noindent[ES10]  A. Edigarian and R.Sigurdsson:  {\it Relative extremal functions and characterization of pluripolar
   sets in complex manifolds}. Trans. Amer. Math. Soc. 362 (2010), no. 10, 5321-5331.

\noindent[KhCh] G. M. Khenkin and E. M. Chirka: {\it Boundary propertie of holomorphic functions of several
 variables}. J. Math. Sci., V.5 (1976), 612-687.

\noindent[Pol99]  E.Poletsky: {\it Disk Envelops of functions}. J. of Functional Analysis V.163 (1999), 111-132.

\noindent [Pin74] S.Pinchuk: {\it A Boundary-uniqueness theorem for holomorphic functions of several complex
variables}. Math. Zametky, V 15, (1974), no.2, 205-212.

\noindent[Sa76] A. Sadullaev: {\it A boundary-uniquiness theorem in $ C^n$}. Mathematical Sbornic, V.101 (143) (1976), no.4, 568-583 = Math.USSR Sb. V.30 (1976), 510-524.

\noindent [Sa80] A. Sadullaev: {\it  P-regularity of sets in $C^n$}. Lect. Not. In Math, V.798 (1980), 402-407.

 \noindent[Sa81] A. Sadullaev: {\it Plurisabharmonic measure and capacity on complex manifolds}. Uspehi Math. Nauk, V.36, (220) (1981), no.34, 53-105= Russian Math. Surveys V.36 (1981), 61-119.

\noindent[Sc62] J. Siciak: {\it On some extremal functions and their
applications in the theory of analytic functions of several
complex variables}. Trans. Amer. Math. Soc. V.105 (1962), 322-357.

%[Si69] J. Siciak: {\it Separatelly analytic functions and envelopes of holomorphy of some
%lower dimensional subsets of $\C^n$}, Ann. Polon. Math. 22 (1969/1970) 145-171

\noindent[Sc81]  J. Siciak: {\it Extremal plurisubharmonic functions in ${\bf C}^{n}$}. Ann. Polon. Math. V.39 (1981), 175-211.

\noindent[Za74] V.P. Zahariuta: {\it Extremal plurisubharmonic
functions,orthogonal polynomials and Bernstein-Walsh theorems for
analytic functions of several complex variables,} I; II. Teor.
Funsttsii Functtsional Anal. i Prilozhen, V.19 (1974), 133-135 ; V.21
(1974), 65-83 (in Russian).

 \noindent[Ze87]. A. Zeriahi: {\it Meilleure approximation polynomiale et croissance des fonctions enti\` eres sur certaines vari\'et\'es alg\'ebriques affines}.  Ann. Inst. Fourier (Grenoble), V.37 (1987), no. 2, 79-104.

\noindent [Ze91]. A. Zeriahi:  {\it Fonction de Green pluricomplexe \` a p\^ole \`a l'infini sur un espace de Stein parabolique et applications.} Math. Scand. V.69 (1991), no. 1, 89-126.\\

\noindent Azimbay Sadullaev\\
 National University of Uzbekistan, \\
100174 Tashkent, Uzbekistan\\
\emph{sadullaev@mail.ru} \\

\noindent Ahmed Zeriahi\\
 Institut de Math\'ematiques de Toulouse \\
Universit\'e Paul Sabatier \\
118 Route de Narbonne, \\
31062 Toulouse \\
\emph{zeriahi@math.univ-toulouse.fr}

\end{document}